% AMSTeX
\def\ver{Sept. 20, 2001, v.5}
\documentstyle{amsppt}
\magnification=1200
\hsize=6.5truein
\vsize=8.9truein
\topmatter
\title Refined Cycle Maps
\endtitle
\author Morihiko Saito
\endauthor
\affil RIMS Kyoto University, Kyoto 606-8502 Japan \endaffil
\abstract
We explain the theory of refined cycle maps associated to arithmetic
mixed sheaves. This includes the case of arithmetic mixed Hodge
structures, and is closely related to work of Asakura, Beilinson,
Bloch, Green, Griffiths, M\"uller-Stach, Murre, Voisin and others.
\endabstract
\endtopmatter
\tolerance=1000
\baselineskip=12pt

\def\scirc{\raise.2ex\hbox{${\scriptstyle\circ}$}}
\def\ssb{\raise.2ex\hbox{${\scriptscriptstyle\bullet}$}}
\def\bl{\bigl}
\def\br{\bigr}
\def\cD{{\Cal D}}
\def\cM{{\Cal M}}
\def\cO{{\Cal O}}
\def\bA{{\Bbb A}}
\def\bC{{\Bbb C}}
\def\bD{{\Bbb D}}
\def\bN{{\Bbb N}}
\def\bP{{\Bbb P}}
\def\bQ{{\Bbb Q}}

\def\bZ{{\Bbb Z}}
\def\bL{\bold{L}}
\def\tM{\widetilde{M}}
\def\tN{\widetilde{N}}
\def\tY{\widetilde{Y}}
\def\ok{\overline{k}}
\def\oM{\overline{M}}
\def\oP{\overline{P}}
\def\oX{\overline{X}}
\def\osig{\overline{\sigma}}
\def\ak{\langle k\rangle}
\def\aok{\langle\overline{k}\rangle}

\def\MM{{\Cal M}{\Cal M}}
\def\SR{\text{\rm SR}}
\def\AJ{\text{\rm AJ}}
\def\go{\text{\rm go}}
\def\ab{\text{\rm ab}}
\def\alg{\text{\rm alg}}
\def\hom{\text{\rm hom}}
\def\dec{\text{\rm dec}}
\def\ind{\text{\rm ind}}
\def\coh{\text{\rm coh}}
\def\hol{\text{\rm hol}}
\def\Hdg{\hbox{{\rm Hdg}}}
\def\Var{\hbox{{\rm Var}}}
\def\LE{\hbox{{\rm LE}}}
\def\Gal{\hbox{{\rm Gal}}}
\def\IC{\hbox{{\rm IC}}}
\def\Perv{\hbox{{\rm Perv}}}
\def\DR{\text{\rm DR}}
\def\CH{\hbox{{\rm CH}}}
\def\Pic{\hbox{{\rm Pic}}}
\def\Spec{\text{{\rm Spec}}\,}
\def\Ext{\hbox{{\rm Ext}}}
\def\Gr{\text{{\rm Gr}}}
\def\Im{\hbox{{\rm Im}}}
\def\Ker{\hbox{{\rm Ker}}}
\def\Coker{\hbox{{\rm Coker}}}
\def\cHom{{\Cal H}om}
\def\Hom{\hbox{{\rm Hom}}}
\def\Sing{\hbox{{\rm Sing}}\,}
\def\supp{\hbox{{\rm supp}}\,}
\def\codim{\hbox{{\rm codim}}}
\def\can{\hbox{{\rm can}}}
\def\div{\hbox{{\rm div}}\,}
\def\MHS{\text{{\rm MHS}}}
\def\MHM{\text{{\rm MHM}}}
\def\an{\text{\rm an}}

\def\SameAuthor{\vrule height3pt depth-2.5pt width1cm}

\document
\centerline{\bf Introduction}

\bigskip\noindent
One of the most fundamental problems in the theory of algebraic cycles
would be Beilinson's conjecture on mixed motives [4], which predicts the
bijectivity of the cycle map
$$
cl : \CH^{p}(X)_{\bQ} \to
{\Ext}_{{D}^{b}\MM(X)}^{2p}(\bQ_{X},\bQ_{X}(p))
\leqno(0.1)
$$
for any smooth projective variety
$ X $ over a field
$ k $.
Here
$ \CH^{p}(X)_{\bQ} $ is the Chow group of codimension
$ p $ algebraic cycles modulo rational equivalence on
$ X $ with
$ \bQ $-coefficients, and
$ D^{b}\MM(X) $ is the bounded derived category of the (conjectural)
abelian category of mixed motivic sheaves on
$ X $.
By the adjoint relation for the structure morphism
$ a_{X} : X \to \Spec k $,
the conjecture would be equivalent to the bijectivity of
$$
cl : \CH^{p}(X)_{\bQ} \to {\Ext}_{{D}^{b}\MM(\Spec k)}^{2p}
(\bQ_{\Spec k}, (a_{X})_{*}\bQ_{X}(p)),
\leqno(0.2)
$$
because
$ \bQ_{X} $ should be the pull-back by
$ a_{X} $ of the constant object
$ \bQ_{\Spec k} $ on
$ \Spec k $.
It is known that this conjecture implies many other important conjectures
on algebraic cycles, such as those of Murre [36], [37], and Bloch [7].

In the case when
$ k $ is embeddable into
$ \bC $ (e.g. if
$ k $ is a number field or
$ \bC) $,
a natural question would be whether
$ \MM(\Spec k) $ is close to the category of
$ \cM_{\SR}(\Spec k) $,
the category of {\it systems of realizations,} which was introduced by
Deligne [20] (see also [21], [30]).
It is expected that the essential image of the natural functor
$ \MM(\Spec k) \to \cM_{\SR}(\Spec
k) $ would be quite close to the full subcategory
$ \cM_{\SR}(\Spec k)^{\go} $ consisting of objects of geometric origin
(which was introduced in [6] for
$ l $-adic sheaves).
So the first test of the conjecture would be whether the cycle map
$$
cl : \CH^{p}(X)_{\bQ} \to {\Ext}_{{D}^{b}{\cM}_{\SR}
(\Spec k)^{\go}}^{2p}(\bQ_{\Spec k},(a_{X})_{*}\bQ_{X}(p))
\leqno(0.3)
$$
is bijective.
We can easily show that the surjectivity of (0.3) is equivalent to the
algebraicity of absolute Hodge cycles [20] on any smooth projective
$ k $-varieties.
See [45], [48].
Thus the surjectivity is essentially reduced to the Hodge conjecture,
which is not easy to prove as well-known.
(However, during an attempt to solve the conjecture we obtained a germ of
a new idea by trying to restrict the Leray spectral sequence to the
generic fiber of a morphism [46].)

Since the category of mixed motives should be universal as far as
cohomology is concerned, we define the category of systems
of realizations as an approximation by endowing the cohomology group with
as much structure as possible.
However, it has been realized by many people that for a complex algebraic
variety (where
$ k = \bC) $,
its cohomology group has more structure.
This was first observed by M. Green, P. Griffiths, and C. Voisin in the
study of the image of the Abel-Jacobi map for a generic hypersurface,
where we have to use the fact that a complex algebraic variety is actually
defined over a finitely generated subring of
$ \bC $.
See [25], [26], [27], [55], [56].
(This fact was also essential for the theory of mod
$ p $ reduction of
$ l $-adic sheaves [6].)
Then, using the models over the finitely generated subrings of
$ \bC $,
it is natural to define the category of {\it arithmetic mixed Hodge
Modules}
(or more generally, {\it arithmetic mixed sheaves}) on a complex algebraic
variety
$ X $ as the inductive limit of the categories of mixed Hodge Modules (or
mixed sheaves) on the models of
$ X $.
In the case the ground field
$ k $ is a finitely generated field over
$ \bQ $, this was already considered in [46] and [47], 1.9
inspired by the arguments in the Appendix to Lect. 1 of [7],
and it is enough to take further the inductive limit over
$ k $.
In particular, we get for the case
$ X = \Spec \bC $ the category of arithmetic mixed Hodge structures,
which is a refinement of mixed Hodge structures
(see also [1] where the theory of mixed sheaves [47] is also used
in an essential way).
It seems that the term ``arithmetic mixed Hodge structure'' was first
used in the work of Green and Griffiths [27], where the
$ \bQ $-structure was not considered because they were mostly interested
in infinitesimal variations of Hodge structures.

The main point of the theory of arithmetic mixed sheaves is that the
injectivity of the refined cycle map of
$ \CH^{2}(X)_{\bQ} $ for a smooth projective complex algebraic variety
$ X $ can be reduced to the injectivity of the Abel-Jacobi map for
codimension two cycles on smooth projective models of
$ X $.
See (4.2).
This shows that an additional hypothesis in [1] is unnecessary.
See also [26], [27], [28].
The result can be extended to some case of higher Chow cycles as in (4.5).
However, it should be noted that this category is too big, and the
forgetful functor to the category of mixed Hodge structures is not fully
faithful (see [51], 2.5).
We would have to restrict to the subcategory consisting of objects of
``geometric origin'' for the study of the surjectivity of the cycle map
although it does not seem to cause a big problem for the injectivity.
Note that our theory applies also to the category of objects of
geometric origin, because it is equivalent to the limit of the category
of objects of geometric origin on the models of
$ X $.
The injectivity of the Abel-Jacobi map for varieties over number fields
has been conjectured by Beilinson [3] and Bloch (see also [9]),
and is one of the most interesting problems in this area.
See also [28].

To show that the obtained new category is really better than the
usual one, we can prove that, restricting the refined cycle map to
the kernel of the Griffiths' Abel-Jacobi map for codimension two cycles,
its image is infinite dimensional if
$ X $ has a nontrivial global two form.
See (4.1).
This was inspired by [58], and is a consequence of Bloch's diagonal
argument in [7] combined with Murre's result on the Albanese motive [36].
(Some special case is treated in [1] using a different method.)
A similar assertion can be proved also for higher Chow groups.
See (4.4).

I would like to thank L. Barbieri-Viale and A. Rosenschon for useful
comments and good questions.
I also thank the referee for useful comments.

We review the theory of mixed Hodge Modules and cycle maps in Sections 1
and 2.
We define the category of arithmetic mixed sheaves in Sect. 3, and states
the main results in Sect. 4.
Some examples are given in Sect. 5.

\bigskip\bigskip
\centerline{{\bf 1. Mixed Hodge Modules}}

\bigskip\noindent
{\bf 1.1.}
$ \cD $-{\bf Modules.} Let
$ X $ be a smooth algebraic variety over a filed
$ k $ of characteristic zero.
Then we have the sheaf of algebraic linear differential operators
$ \cD_{X/k} $ which has the increasing filtration
$ F $ by the order of operators.
Since
$ \cD_{X/k} = \cD_{X/k'} $ for a finite extension
$ k' \subset k $,
we will write
$ \cD_{X} $ for
$ \cD_{X/k} $ in the sequel.

Let
$ (M,F) $ be a coherent filtered
$ \cD_{X} $-Module.
Here we assume always that
$ F $ is exhaustive and
$ F_{p}M = 0 $ for
$ p \ll 0 $.
Then the coherence of
$ (M,F) $ means that
$ \Gr^{F}M \,(:= \oplus_{p} \Gr_{p}^{F}M) $ is coherent over
$ \Gr^{F}\cD_{X} $.
We say that
$ (M,F) $ is {\it holonomic} if it is coherent and
$ \dim \supp \Gr^{F}M = \dim X $,
where
$ \supp \Gr^{F}M $ is a subvariety of the cotangent space of
$ X $ which is naturally isomorphic to
$ \Spec \Gr^{F}\cD_{X} $.
We say that a
$ \cD_{X} $-Module is holonomic if it is a coherent
$ \cD_{X} $-Module and having locally a filtration
$ F $ such that
$ (M,F) $ is holonomic.
The category of coherent (or holonomic) filtered
$ \cD_{X} $-Modules will be denoted by
$ MF_{\coh}(\cD_{X}) $ (or
$ MF_{\hol}(\cD_{X})) $.
Forgetting the filtration, we have
$ M_{\coh}(\cD_{X}), M_{\hol}(\cD_{X}) $ similarly.

If
$ k = \bC $ and
$ A $ is a subfield of
$ \bC $,
let
$ \Perv(X,A) $ denote the abelian category of perverse sheaves on
$ X^{\an} $ with
$ A $-coefficients and with algebraic stratifications.
See [6].
Then we have the de Rham functor
$$
\DR : M_{\hol}(\cD_{X}) \to \Perv(X,\bC)
\leqno(1.1.1)
$$
defined by
$ M \to \Omega_{X}^{\dim X}{\otimes}_{{\cD}_{X}}^{\bL}M $ (using a
natural projective resolution of the right
$ \cD_{X} $-Module
$ \Omega_{X}^{\dim X}) $.
See e.g. [11].
This functor is exact and faithful.

\medskip\noindent
{\bf 1.2.~Direct images.} Let
$ f : X \to Y $ be a proper morphism of smooth
$ k $-varieties.
Then the cohomological direct image functor
$ H^{j}f_{*} : MF_{\coh}(\cD_{X}) \to MF_{\coh}(\cD_{Y}) $ is defined
by factorizing
$ f $ by
$ X \to X \times_{k} Y \to Y $,
where the first morphism is the embedding by the graph of
$ f $,
and the second is the projection.

If
$ i : X \to Y $ is a closed embedding of codimension
$ d $,
let
$ (y_{1}, \dots, y_{n}) $ be a local coordinate system of
$ Y $ (i.e. it defines an \'etale map of an open subvariety of
$ Y $ to
$ \bA^{n}) $ such that
$ X = \{y_{1} = \dots = y_{d} = 0\} $.
Let
$ \partial_{j} = \partial /\partial y_{j} $ so that
$ \cD_{Y} $ is locally identified with
$ \cO_{Y}\otimes_{\bC}\bC[\partial_{1},
\dots, \partial_{n}] $.
Then the direct image
$ i_{*}M $ is locally isomorphic to
$ M\otimes_{\bC}\bC[\partial_{1}, \dots, \partial_{d}]
 $ using the coordinates (see [11]), and the filtration
$ F $ on
$ i_{*}M $ is defined by
$$
F_{p}i_{*}M = \sum_{\nu \in \bN^{d}} F_{p-|\nu |-d}M\otimes
\partial^{\nu},
\leqno(1.2.1)
$$
where
$ \partial^{\nu} = \prod_{1\le k\le d} {\partial}_{k}^{{\nu}_{k}} $ for
$ \nu = (\nu_{1}, \dots, \nu_{d}) \in \bN^{d} $.
(This is independent of the choice of the coordinates.)

If
$ q : X \times_{k} Y \to Y $ is the projection, we have the relative
de Rham complex
$ \DR_{X\times Y/Y}(M,F) $ such that
$ F_{p}\DR_{X\times Y/Y}(M)^{j} = \Omega_{X\times Y/Y}^{j}\otimes_{\cO}
F_{p+j}M $.
Then the cohomological direct image is defined to be the cohomology sheaf
of the filtered direct image of
$ \DR_{X\times Y/Y}(M,F) $ by
$ q $.
(Note that the filtration on the direct image complex is not necessarily
strict, and we take the induced filtration on the cohomology sheaf.)

It is known that holonomic filtered
$ \cD $-Modules are stable by the cohomological direct image under
a proper morphism, and the cohomological direct image functors for
(filtered)
$ \cD $-Modules and perverse sheaves are compatible with each other
via the functor (1.1.1), see e.g. [11], [43].

\medskip\noindent
{\bf 1.3.~Mixed Hodge Modules on complex algebraic varieties.} Let
$ X $ be a smooth complex algebraic variety.
We say that
$ M = ((M_{\cD},F) $,
$ M_{\bQ} $,
$ W $,
$ \alpha ) $ is a bifiltered holonomic
$ \cD_{X} $-Module with rational structure, if
$ (M_{\cD},F) $ is a holonomic filtered
$ \cD_{X} $-Module,
$ M_{\bQ} $ is a perverse sheaf with rational coefficients,
$ \alpha $ is a comparison isomorphism
$ \DR(M_{\cD}) = M_{\bQ}\otimes_{\bQ}\bC
 $ in
$ \Perv(X,\bC) $,
and
$ W $ is a pair of finite increasing filtrations on
$ M_{\cD} $ and
$ M_{\bQ} $ compatible with
$ \alpha $.
Here
$ \Perv(X,A) $ for
$ A = \bQ $ or
$ \bC $ is the full subcategory of
$ \Perv(X^{\an},A) $ (see [6]) consisting of perverse sheaves with
algebraic stratifications, and similarly for
$ {D}_{c}^{b}(X,A) $.
A mixed Hodge Module on
$ X $ is a bifiltered holonomic
$ \cD_{X} $-Module
$ M $ with rational structure satisfying the following conditions:

The first condition is that mixed Hodge Modules are defined
Zariski-locally.
The second is that, restricting a mixed Hodge Module
$ M $ to an open smooth subvariety
$ Z $ of
$ \supp M $ on which
$ K $ is a (shifted) local system, it is isomorphic to the direct image by
the closed immersion
$ Z \to X \setminus (\supp M \setminus Z) $ of an admissible
variation of mixed Hodge structure in the sense of [32], [54] (see also
(1.4) below) up to a shift of the weight filtration.
Here the converse is also true, and an admissible variation of mixed Hodge
structure is a mixed Hodge Module ([44], 3.27).
The last condition claims that mixed Hodge Modules are locally obtained
inductively by gluing mixed Hodge Modules supported on a divisor and
admissible variations of mixed Hodge structures on a smooth closed
subvariety in the complement on the divisor:

Let
$ g $ be a function on
$ X $ such that the restriction of a bifiltered holonomic
$ \cD_{X} $-Module with rational structure
$ M $ to the complement
$ U $ of
$ g^{-1}(0) $ is the direct image of an admissible variation of mixed
Hodge structure
$ M' $ on a smooth variety
$ Y $ by the closed immersion
$ i : Y \to U $.
Let
$ \psi_{g,1} $ and
$ \varphi_{g,1} $ be the nearby and vanishing cycle functors with
unipotent monodromy [18].
Then the condition is that the nearby and vanishing cycles
$ \psi_{g,1}M, \varphi_{g,1}M $ are well defined and the obtained
$ M'' := \varphi_{g,1}M $ is a mixed Hodge Module supported on
$ g^{-1}(0) $.
(The first condition consists of the compatibility of the three
filtrations
$ F, W, V $ and the existence of the relative monodromy filtration,
see [44], 2.2.)

Here we have canonical morphisms of mixed Hodge Modules
$$
\can : \psi_{g,1}i_{*}M' \to M'',\quad
\Var : M'' \to \psi_{g,1}i_{*}M'(-1),
$$
satisfying the gluing condition
$$
\Var\scirc\can = N,
\leqno(1.3.1)
$$
where
$ N = \log T_{u} $ with
$ T = T_{s}T_{u} $ the Jordan decomposition of the monodromy
$ T $.
(Here
$ (-1) $ denotes the Tate twist.)
We can show that
$ M $ is uniquely determined by
$ (M', M'', \can, \Var) $ satisfying the gluing condition (1.3.1).
More precisely, we have an equivalence between the category of mixed Hodge
Modules on
$ X $ and the category of
$ (M', M'', \can, \Var) $ satisfying (1.3.1).
See [44], 2.28.
Furthermore, the corresponding
$ M $ is uniquely determined by using the above condition on the
well-definedness of the nearby and vanishing cycle functors.
See [44], 2.8.

We denote by
$ \MHM(X) $ the category of mixed Hodge Modules on
$ X $.
It is an abelian category such that any morphism is strictly compatible
with the Hodge and weight filtrations
$ F, W $ in the strong sense.
See [43], 5.1.14.

Now let
$ X $ be a (singular) complex algebraic variety.
We consider the category
$ \LE(X) $ whose objects are closed embeddings
$ U \to V $ where
$ U $ is an open subvariety of
$ X $ and
$ V $ is smooth.
The morphisms are pairs of morphisms between
$ U $ and between
$ V $ making a commutative diagram.
Here the morphisms of
$ U $ are assumed to be compatible with the inclusions to
$ X $.
For
$ \{U \to V\} \in \LE(X) $,
let
$ \MHM_{U}(V) $ denote the category of mixed Hodge Modules on
$ V $ supported on
$ U $.
Then a mixed Hodge Module on
$ X $ is a collection of mixed Hodge Modules
$ M_{U\to V} \in \MHM_{U}(V) $ for
$ \{U \to V\} \in \LE(X) $ (which is called the representative of
$ M $ for
$ \{U \to V\}) $ together with isomorphisms
$$
v_{*}M_{U\to V}|_{V'\setminus (U'\setminus U)} =
M_{U'\to V'}|_{V'\setminus (U'\setminus U)}
$$
for
$ (u,v) : \{U \to V\} \to \{U' \to V'\} $ satisfying the
usual cocycle condition.
(We can define the category of filtered
$ \cD $-Modules on
$ X $ similarly.)

In the case
$ X $ is smooth, we can show that this definition is equivalent to the
previous one (i.e. we have naturally an equivalence of categories).

Actually, to define a mixed Hodge Module on a singular
$ X $,
it is not necessary to define
$ M_{U\to V} $ for all
$ \{U \to V\} $;
it is enough to do so for
$ \{U \to V\} $ such that the
$ U $ cover
$ X $,
but the gluing morphisms are defined by using the closed embeddings
$ U \cap U' \to V \times V' $.

\medskip\noindent
{\bf 1.4.~Admissible variation of mixed Hodge structure.} Let
$ X $ be a smooth complex algebraic variety, and
$ \oX $ a smooth compactification of
$ X $ such that the complement
$ D := \oX \setminus X $ is a divisor with normal crossings.
Let
$ M = ((M_{\cO};F,W),(M_{\bQ},W)) $ be a variation of mixed Hodge
structure on
$ X $,
where
$ M_{\cO} $ is the underlying
$ \cO_{X} $-Module with the integral connection
$ \nabla $ having regular singularity at infinity, and
$ M_{\bQ} $ is the underlying
$ \bQ $-local system with an isomorphism
$ \Ker\, \nabla^{\an} = M_{\bQ}\otimes_{\bQ}\bC
 $.
By [32], [54],
$ M $ is admissible if the graded pieces
$ \Gr_{k}^{W}M $ are polarizable variations of Hodge structures
with quasi-unipotent local monodromies, and
furthermore the following conditions are satisfied:

In the case the local monodromies around
$ D $ are all unipotent, let
$ \oM_{\cO} $ be Deligne's extension of
$ M_{\cO} $ (see [17]).
Then

\noindent
(i) The filtrations
$ F, W $ on
$ M $ are extended to
$ \oM_{\cO} $ so that
$ \Gr_{F}^{p}\Gr_{k}^{W}\oM_{\cO} $ is a locally free
$ \cO_{\oX} $-Module for any
$ p, k $.

\noindent
(ii) The relative monodromy filtration exists for the local monodromy
around each irreducible component of
$ D $.

Here we can replace
$ \oM_{\cO} $ with its restriction to
$ \oX \setminus \Sing D $,
i.e. it is enough to consider the conditions around the smooth points of
$ D $,
see [32].

In general, the above conditions should be satisfied for the pull-back of
$ M $ by a dominant morphism such that the local monodromies of the
pull-back are unipotent.
More precisely, the condition for admissible variation is analytic local
on a compactification of
$ X $,
and it is enough to take locally a ramified cover of a polydisk as usual.
Here we can consider
$ {M}_{\cO}^{\an} $ instead of
$ M_{\cO} $,
because
$ M_{\cO} $ is uniquely determined by
$ {M}_{\cO}^{\an} $ due to the regularity and GAGA.

\medskip\noindent
{\bf 1.5.~Mixed Hodge Modules on algebraic varieties.} Let
$ k $ be a subfield of
$ \bC $,
and
$ X $ a smooth
$ k $-variety.
Let
$ X_{\bC} = X\otimes_{k}\bC $.
Then a mixed Hodge Module
$ M $ on
$ X $ consists of
$ ((M_{\cD},F),W), (M_{\bQ},W) $ and
$ \alpha $ such that
$ (((M_{\cO}\otimes_{k}\bC,F),W), (M_{\bQ},W), \alpha ) $ is a mixed
Hodge Module on
$ X_{\bC} $,
where
$ (M_{\cO},F) \in MF_{\hol}(\cD_{X}), M_{\bQ} \in
\Perv(X_{\bC},\bQ) $ with a finite increasing filtration
$ W $,
and
$ \alpha $ is a comparison isomorphism
$ \DR(M\otimes_{k}\bC) = K\otimes_{\bQ}\bC
 $ in
$ \Perv(X_{\bC},\bC) $ which is compatible with
$ W $.
Here we also assume that polarizations on the graded pieces
$ \Gr_{k}^{W}M $ are defined over
$ k $,
i.e., they are induced by isomorphisms of filtered
$ \cD_{X} $-Modules
$ \Gr_{k}^{W}(M_{\cD},F)(k) \simeq \bD\Gr_{k}^{W}(M_{\cD},F) $
compatible with pairings of perverse sheaves, where
$ \bD $ denotes the dual.
This is necessary to assure that the graded pieces are semisimple.

We will denote by
$ \MHM(X/k) $ the category of mixed Hodge Modules on
$ X/k $.
This is an abelian category such that every morphism is strictly
compatible with
$ F, W $ in the strong sense.
We have naturally the forgetful functors
$$
\MHM(X/k) \to \MHM(X_{\bC}) \to \Perv(X_{\bC},\bQ),
$$
which induce
$$
D^{b}\MHM(X/k) \to D^{b}\MHM(X_{\bC}) \to
{D}_{c}^{b}(X_{\bC},\bQ),
\leqno(1.5.1)
$$
using the canonical functor
$ D^{b}\Perv(X_{\bC},\bQ) \to {D}_{c}^{b}(X_{\bC},\bQ) $ in [6].

We can define similarly the notion of admissible variation of mixed Hodge
structures on
$ X $ (i.e., the underlying filtered
$ \cO $-Modules and polarizations are defined over
$ X $.)
Then we see that mixed Hodge Modules on
$ X $ are obtained locally by gluing mixed Hodge Modules supported on a
divisor and admissible variations of mixed Hodge structure on a smooth
closed subvariety in the complement of the divisor as in (1.3).

\medskip\noindent
{\bf 1.6.~Mixed sheaves.} In this paper we consider more generally the
category of mixed sheaves
$ \cM(X/k) $ in the sense of [47].
However, to simplify the explanation, we assume in this paper that
$ \cM(X/k) $ is either
$ \MHM(X/k) $ defined above or the category
$ \cM_{\SR}(X/k) $ consisting of {\it systems of realizations}
$ ((M_{\cD};F,W), (M_{\sigma},W), (M_{l},W)) $,
where
$ (M_{\cD},F) $ is a holonomic filtered
$ \cD $-Module on
$ X $ endowed with a finite filtration
$ W $,
$ (M_{\sigma},W) $ for an embedding
$ \sigma : k \to \bC $ is a filtered perverse sheaf on
$ (X\otimes_{k,\sigma}\bC)^{\an} $ with
$ \bQ $-coefficients, and
$ (M_{l},W) $ for a prime number
$ l $ is a filtered perverse
$ l $-adic sheaf on
$ X_{\ok} := X\otimes_{k}\ok $ with
$ \bQ_{l} $-coefficients which has a continuous action of the
Galois group of
$ \ok/k $ (i.e. the action is lifted to perverse sheaves with
$ \bZ_{l} $-coefficients).
Furthermore these are endowed with comparison isomorphisms
$$
\DR((M_{\cD},W)\otimes_{k,\sigma}\bC) = (M_{\sigma},
W)\otimes_{\bQ}\bC,\quad \epsilon^{*}{i}_{\osig}^{*}
(M_{l},W) = (M_{\sigma},W)\otimes_{\bQ}\bQ_{l}
$$
for an extension
$ \osig : \ok \to \bC $ of
$ \sigma $ (in a compatible way with the action of
$ \Gal(\ok/k), $ see [30]).
Here
$ A = \bQ, \ok $ is the algebraic closure of
$ k $ in
$ \bC $,
and
$ i_{\osig} $:
$ X_{\bC} \to X_{\ok} $ is the canonical morphism.
(See [6] for
$ \epsilon^{*} $.)
In the case
$ X = \Spec k $,
$ \cM_{\SR}(\Spec k) $ coincides with the category of systems of
realizations introduced by Deligne [20], [21] (and this formulation
is due to Jannsen [30]).

For
$ \cM(X/k) = \MHM(X/k) $ or
$ \cM_{\SR}(X/k) $,
there exists canonically the base change functor
$$
\cM(X/k) \to \cM(X\otimes_{k}k'/k')
\leqno(1.6.1)
$$
for a finite extension
$ k \subset k' $ in a compatible way with the cohomological direct image
and pull-back and also with dual and external product, etc.
There is also the (canonically defined) forgetful functor
$$
\cM(X/k) \to \MHM(X/k),
\leqno(1.6.2)
$$
which is compatible with the standard functors as above.

If
$ X/k $ is smooth and purely
$ d $-dimensional, we have the constant object
$ \bQ_{X/k} $ in
$ D^{b}\cM(X/k) $ which actually belongs to
$ \cM(X/k)[-d] $ by the definition of the perverse sheaf [6].
We have also the Tate twist
$ \bQ_{X/k}(j) $ for
$ j \in \bZ $ (using the cohomology of
$ \bP^{-j} $ for
$ j $ negative, and taking the dual for
$ j $ positive).

\medskip\noindent
{\bf 1.7.~Theorem.} {\it
For a morphism
$ f $ of algebraic varieties over
$ k $,
there exist canonical functors
$ f_{*} $,
$ f_{!} $,
$ f^{*} $,
$ f^{!} $,
$ \bD $,
$ \otimes $,
$ \cHom $,
etc. between the derived categories
$ D^{b}\cM(X/k) $ in a compatible way with the corresponding functors
between the derived categories
$ D^{b}\MHM(X/k), D^{b}\MHM(X_{\bC}) $ or
$ {D}_{c}^{b}(X_{\bC},\bQ) $ via the functors (1.5.1) and (1.6.2).
}

\medskip\noindent
{\it Proof.}
This follows from the same argument as in [44].
See also [47].
Indeed, using Beilinson's resolution, the stability by the direct image
is reduced to the one by the cohomological direct image for an affine
morphism [5].
(If
$ f $ is quasi-projective, this is especially simple by taking two sets
of affine coverings of
$ X $ associated with general hyperplane sections and using the co-Cech
and Cech complexes.)
The pull-backs are defined to be the adjoint functors of the
direct images.
For the existence, we may assume
$ f $ is either a closed embedding
$ i $ or a projection
$ p $.
In the former case, the assertion is reduced to the full faithfulness of
the direct image
$$
i_{*} : D^{b}\cM(X/k,\bQ) \to D^{b}\cM(Y
/k,\bQ),
$$
which is shown by using the functor
$ \xi_{g} $ in [44], 2.22.
In the latter case it is enough to show the existence of
$ {a}_{X}^{*}\bQ_{\Spec k/k} $ (using the duality and the external
product).
But this is represented by any complex
$ M^{\ssb} $ having a morphism
$ \bQ_{\Spec k/k} \to (a_{X/k})_{*}M^{\ssb} $ in
$ D^{b}\cM(X/k) $ such that the image of
$ M^{\ssb} $ in
$ {D}_{c}^{b}(X_{\bC},\bQ) $ is isomorphic to
$ \bQ_{X_{\bC}} $ and the image of the morphism is identified
with the canonical morphism
$ \bQ \to (a_{X_{\bC}})_{*}\bQ_{X_{\bC}} $.
So it exists locally on open subsets which are embeddable into smooth
varieties, and we can glue them by using the adjoint morphism for the
inclusion of open subvarieties, see [44], 4.4.

\medskip\noindent
{\bf 1.8.~Definition.} For a
$ k $-variety
$ X $ with structure morphism
$ a_{X/k} : X \to \Spec k $,
we define
$$
\bQ_{X/k}(j) = {a}_{X/k}^{*}\bQ_{\Spec k/k}(j),\quad
H^{j}(X/k,\bQ(j))=H^{j}(a_{X/k})_{*}\bQ_{X/k}(j).
\leqno(1.8.1)
$$
We omit
$ /k $ in the case
$ k = \bC $.

\medskip\noindent
{\bf 1.9.~Decomposition of the direct images.} If
$ X $ is smooth proper over
$ k $,
we have a noncanonical isomorphism
$$
(a_{X/k})_{*}\bQ_{X/k} \simeq \oplus_{j} H^{j}(X/k,\bQ)[-j]\quad
\text{in }D^{b}\cM(\Spec k/k).
\leqno(1.9.1)
$$
See e.g. [44], 4.5.3.

\bigskip\bigskip
\centerline{{\bf 2. Cycle Map and Geometric Origin}}

\bigskip\noindent
{\bf 2.1.~Cycle map.} Let
$ X $ be a smooth
$ k $-variety.
We define an analogue of Deligne cohomology by
$$
\aligned
{H}_{\cD}^{i}(X/k,\bQ(j))
&= \Ext^{i}(\bQ_{X/k}, \bQ_{X/k}(j))
\\
\bl(
&= \Ext^{i}(\bQ_{k}, (a_{X/k})_{*}\bQ_{X/k}(j))\br),
\endaligned
\leqno(2.1.1)
$$
where the extension groups are taken in the derived category of
$ \cM(X/k) $ or
$ \cM(\Spec k/k) $,
and the second isomorphism follows from the adjoint relation between the
direct image and the pull-back by
$ a_{X/k} $.
(In the case
$ k = \bC $,
we will often omit
$ /k $ to simplify the notation.)

Let
$ \CH^{p}(X)_{\bQ} $ be the Chow group consisting of codimension
$ p $ cycles modulo rational equivalence on
$ X $ with rational coefficients.
Then we have naturally the cycle map
$$
cl : \CH^{p}(X)_{\bQ} \to {H}_{\cD}^{2p}(X/k,\bQ(p)).
\leqno(2.1.2)
$$
If a cycle
$ \zeta $ is represented by an irreducible closed subvariety
$ Z $,
then
$ cl(\zeta ) $ is defined to be the composition of
$$
\bQ_{X/k} \to \bQ_{Z/k} \to \IC_{Z/k}\bQ
[-d_{Z/k}]
$$
with its dual, by using the dualities
$$
\bD(\bQ_{X/k}) = \bQ_{X/k}(d_{X/k})[2d_{X/k}],\quad \bD
(\IC_{Z/k}\bQ) = \IC_{Z/k}\bQ(d_{Z/k}),
$$
where
$ \IC_{Z/k}\bQ $ is the intersection complex, and
$ d_{X/k} = \dim X/k $.
See [44], 4.5.15.
We can show that the cycle map is compatible with the pushdown and the
pull-back of cycles (where a morphism is assumed to be proper in the case
of pushdown.)
See [45, II].

The composition of (2.1.2) with the natural projection
$$
\CH^{p}(X)_{\bQ} \to \Hom\bl(\bQ_{k}, H^{2p}(X/k,\bQ(p))\br)
$$
is the usual cycle map.
Let
$ \CH_{\hom}^{p}(X)_{\bQ} $ be its kernel (which consists of
homologically equivalent to zero cycles).
Then (2.1.2) induces a generalized Abel-Jacobi map over
$ k $:
$$
\CH_{\hom}^{p}(X)_{\bQ} \to J^{p}(X/k)_{\bQ} :=
\Ext^{1}\bl(\bQ_{k}, H^{2p-1}(X/k,\bQ(p))\br),
\leqno(2.1.3)
$$
where
$ \Ext^{1} $ is taken in
$ \cM(\Spec k/k) $.

Let
$ \CH^{p}(X,m)_{\bQ} $ be Bloch's higher Chow group with rational
coefficients of a smooth
$ k $-variety
$ X $ [10].
By [47] we have the cycle map
$$
cl : \CH^{p}(X,m)_{\bQ} \to {H}_{\cD}^{2p-m}(X/k,\bQ(p)).
\leqno(2.1.4)
$$
If
$ X $ is smooth proper over
$ k $ and
$ m > 0 $,
then this cycle map induces a generalized Abel-Jacobi map over
$ k $:
$$
\CH^{p}(X,m)_{\bQ} \to \Ext^{1}\bl(\bQ_{k},
H^{2p-m-1}(X/k, \bQ(p))\br),
\leqno(2.1.5)
$$
because the Leray spectral sequence for
$ {H}_{\cD}^{2p-m}(X/k,\bQ(p)) $ degenerates at
$ E_{2} $ by (1.9.1), and
$ \Hom\bl(\bQ_{k}, H^{2p-m}(X/k, \bQ(p))\br) = 0 $.

\medskip\noindent
{\bf 2.2.~Griffiths' Abel-Jacobi map.} If
$ k = \bC $ and
$ \cM(X) = \MHM(X) $,
then
$ {H}_{\cD}^{i}(X,\bQ(j)) $ for a smooth projective variety
$ X $ coincides with Deligne cohomology in the usual sense ([22],
$ [23]), $ and the cycle map (2.1.2) coincides with Deligne's cycle map.
In particular, (2.1.3) coincides with Griffiths' Abel-Jacobi map
$$
\CH_{\hom}^{p}(X) \to J^{p}(X)
\,\bl(= {\Ext}_{\MHS}^{1}\bl(\bZ, H^{2p-1}(X,\bZ(p))\br)\br)
\leqno(2.2.1)
$$
tensored with
$ \bQ $, where
$ J^{p}(X) $ is the Griffiths intermediate Jacobian [29],
and the last isomorphism follows from [12].

\medskip\noindent
{\bf 2.3.~Injectivity of the Abel-Jacobi map.}
It is expected that higher extension groups
$ \Ext^{i} $
$ (i > 1) $ should vanish in the (conjectural) category of mixed motives
over a number field.
Since the category of systems of realizations is an approximation of the
category of mixed motives, it is interesting whether the Abel-Jacobi map
(2.1.3) is injective in the case
$ k $ is a number field.
Actually Beilinson conjectures the injectivity of the composition of
(2.1.3) with the natural morphism:
$$
\CH_{\hom}^{p}(X)_{\bQ} \to J^{p}(X/k)_{\bQ} \to J^{p}(X_{\bC})_{\bQ}
\leqno(2.3.1)
$$
at least if we restrict it to the subgroup
$ \CH_{\alg}^{p}(X)_{\bQ} $ consisting of algebraically equivalent to
zero cycles [3].
Since the image of
$ \CH_{\alg}^{p}(X)_{\bQ} $ by (2.1.3) is contained in the algebraic
part of the Jacobian, and the restriction of (2.3.1) to
$ \CH_{\alg}^{p}(X)_{\bQ} $ is defined algebraically (3.10),
it would be natural to conjecture the injectivity of (2.3.1) for the
algebraically equivalent to zero cycles.
However, it may be better to conjecture the injectivity of (2.1.3)
in general.

\medskip\noindent
{\bf 2.4.~Geometric origin.} We denote by
$ \cM(X/k)^{\go} $ the full subcategory of
$ \cM(X/k) $ consisting of objects of geometric origin (see [6] for the
case of perverse sheaves).
This is by definition the smallest full subcategory of
$ \cM(X/k) $ which is stable by the standard cohomological functors
$ H^{j}f_{*}, H^{j}f_{!}, H^{j}f^{*}, H^{j}f^{!}, $ etc. and also by
subquotients in
$ \cM(X/k) $,
and contains the constant object
$ \bQ_{\Spec k/k} $ for
$ X = \Spec k $.
(This satisfies the axiom of mixed sheaves, see [47], 7.1.)
Actually, it is enough to assume the stability by the cohomological direct
images and pull-backs, because the nearby and vanishing cycle functors are
expressed by using the direct images and pull-backs in the same way as in
[18] (see e.g. [47], 5.7), and the stability by dual and external
product follows from the compatibility with those functors.
(More precisely we have (2.5) below.)

We define
$ {H}_{\cD}^{i}(X/k,\bQ(j))^{\go} $ as in (2.1.1) with
$ \cM(X/k) $ replaced by
$ \cM(X/k)^{\go} $.
Note that the natural morphism
$$
{H}_{\cD}^{i}(X/k,\bQ(j))^{\go} \to {H}_{\cD}^{i}(X/k,\bQ(j))
\leqno(2.4.1)
$$
is not injective in general.
We have the cycle map
$$
cl : \CH^{p}(X)_{\bQ} \to {H}_{\cD}^{2p}(X/k,\bQ(p))^{\go}
\leqno(2.4.2)
$$
factorizing (2.1.2).
Let
$$
\aligned
J^{p}(X/k)_{\bQ}
&= {\Ext}_{\cM}^{1}\bl(\bQ, H^{2p-1}(X/k,\bQ(p))\br),
\\
J^{p}(X/k)_{\bQ}^{\go}
&= {\Ext}_{{\cM}^{\go}}^{1}\bl(\bQ, H^{2p-1}(X/k,\bQ(p))\br),
\endaligned
$$
where
$ \cM $ and
$ \cM^{\go} $ mean
$ \cM(\Spec k/k) $ and
$ \cM(\Spec k/k)^{\go} $.
We have a canonical injection
$$
J^{p}(X/k)_{\bQ}^{\go} \to J^{p}(X/k)_{\bQ}.
\leqno(2.4.3)
$$
(If
$ k = \bC $ and
$ \cM(X) = \MHM(X) $,
we omit
$ /k $.)

\medskip

We can show the following:

\medskip\noindent
{\bf 2.5.~Proposition.} {\it
For
$ \cM \in \cM(X/k) $,
it is of geometric origin if and only if for any point of
$ X $,
there exist an open neighborhood
$ U $,
a closed embedding
$ i : U \to Z $,
a quasi-projective morphism
$ \pi : Y \to Z $,
and a divisor
$ D $ on
$ Y $ such that
$ i_{*}\cM|_{U} $ is isomorphic to a subquotient of
$ H^{j}\pi_{*}j_{!}\bQ_{(Y\setminus D)/k} $ in
$ \cM(Z/k) $.
Here
$ j : Y \setminus D \to Y $ denotes the inclusion morphism, and
we may assume that
$ D $ is a divisor with normal crossings on
$ Y $.
(See [47] and also [45, I].)
}

\medskip\noindent
{\bf 2.6.~Theorem.} {\it
Assume
$ \cM(X/k) = \cM_{\SR}(X/k) $ in (1.6).
Then the following assertions are equivalent:

\noindent
(i) The cycle map (2.4.2) is surjective for any smooth projective variety
$ X $ over
$ k $.

\noindent
(ii) Absolute Hodge cycles on any smooth projective varieties over
$ k $ are algebraic.

If
$ k = \bC $ and
$ \cM(X) = \MHM(X) $,
then the equivalence holds with absolute Hodge cycles replaced by Hodge
cycles, and these are further equivalent to:

\noindent
(iii) The images of (2.2.1) and (2.4.3) coincide for any smooth complex
projective varieties.
}

\medskip
(See [45, I] and [48].)

\medskip\noindent
{\bf 2.7.~Remark.} By (2.6), the surjectivity of the cycle map (2.4) is
reduced to the algebraicity of absolute Hodge cycles, and the latter is
easily reduced to the usual Hodge conjecture.
To show the last conjecture, it would be natural to consider a morphism
of a variety (which is birational to
$ X) $ to another variety of dimension
$ \ge 2 $ (e.g. by taking Lefschetz pencils successively).
In [46] we tried to restrict the Leray spectral sequence to each fiber
using extension groups.
Here the category of usual mixed Hodge structures is not good enough
because of the vanishing of higher extension groups.
But it is also unclear whether the full subcategory of objects of
geometric origin is useful, since their higher extension groups
are very difficult to calculate, although they are not expected
to vanish.

Another big problem in this attempt is that, even if we could get cycles
on fibers, it is not clear whether they come from one cycle on the total
space.
The difficulty comes from the fact that the higher extension classes do
not form a section of some geometric object over the base space as in the
case of normal functions.

To solve the last problem it would be natural to consider the generic
fiber of the morphism and try to find an argument corresponding to
``spreading out'' of algebraic cycles (as in [7], p. 1.20).
Then we get the idea of taking the inductive limit of the category of
mixed sheaves on the pull-back of nonempty open subvarieties of the base
space.
See [46] and [47], 1.9.
Unfortunately, this idea did not work well for the problem mentioned
above, because we do not yet have a good category of mixed sheaves (which
is strong enough to solve the problem).
However, seeing earlier work of M. Green [25] and C. Voisin [55], we
notice that the cohomology of a complex algebraic variety has more
structure than expressed in the systems of realizations, and find that
the converse of the above argument would be possible.
A complex algebraic variety
$ X $ has a model of finite type over a subfield
$ k $ of
$ \bC $ having a morphism to another integral
$ k $-variety whose geometric generic fiber over
$ \Spec \bC $ is isomorphic to
$ X $.
Then it is natural to consider the inductive limit of mixed sheaves on the
models of
$ X $.
See also [1] and especially [56], p. 194.
(Note that a similar idea was essentially used in the theory of mod
$ p $ reduction for perverse sheaves [6].)
Thus we get the notion of arithmetic mixed sheaf which will be explained
in the next section.

\bigskip\bigskip
\centerline{{\bf 3. Arithmetic Mixed Sheaves}}

\bigskip\noindent
{\bf 3.1.~Construction.} Let
$ k, K $ be subfields of
$ \bC $ such that
$ k \subset K $.
Then for a
$ K $-variety
$ X $,
there exists a finitely generated
$ k $-subalgebra
$ R $ of
$ K $ such that
$ X $ is defined over
$ R $,
i.e., there is an
$ R $-scheme
$ X_{R} $ of finite type such that
$ X = X_{R}\otimes_{R}K $.
For a finitely generated smooth
$ k $-subalgebra
$ R' $ of
$ K $ containing
$ R $,
let
$ X_{R'} = X_{R}\otimes_{R}R', S' = \Spec R', d_{R'} = \dim_{k} S' $,
and let
$ k_{R'} $ be the algebraic closure of
$ k $ in
$ R' $.
Then
$ k_{R'} $ coincides with the algebraic closure of
$ k $ in the function field
$ k(S') $ of
$ S' $ (because
$ S' $ is normal), and
$ k_{R'} $ is a finite extension of
$ k_{R} $.
Furthermore,
$ S/k_{R'} $ is geometrically irreducible and
$ S'_{\bC} := S'\otimes_{k_{R'}}\bC $ is connected.
We define
$$
\cM(X/K)_{\ak} = \varinjlim \cM(X_{R'}/k_{R'})[-d_{R'}],
$$
where the inductive limit is taken over
$ R' $ as above.
Here we denote by
$ \cM(X_{R'}/k_{R'})[-d_{R'}] $ the category of mixed sheaves
$ \cM(X_{R'}/k_{R'}) $ shifted by
$ -d_{R'} $ in the derived category.
Note that the shift is necessary due to the normalization of perverse
sheaves in [6] which implies that perverse sheaves are stable by the
usual pull-back
$ f^{*} $ under a smooth morphism
$ f $ up to the shift of complexes by the relative dimension.

More precisely, we define the order relation
$ R' < R'' $ by the inclusion
$ R' \subset R'' $ together with the smoothness of
$ R'' $ over
$ R' $.
Then we have natural functors
$$
\aligned
\cM(X_{R'}/k_{R'})[-d_{R'}]
&\to \cM(X_{R'}\otimes_{k_{R'}}k_{R''}/k_{R''})[-d_{R'}]
\\
&\to \cM(X_{R''}/k_{R''})[-d_{R''}],
\endaligned
\leqno(3.1.1)
$$
where the first comes from (1.6.1).
Note that
$ R'\otimes_{k_{R'}}k_{R''} \to R'' $ is injective because
$ k_{R'} $ is algebraically closed in the fraction field of
$ R'. $

We have the canonical functors
$$
\cM(X_{R'}/k_{R'})[-d_{R'}] \to \MHM(X_{\bC})
\leqno(3.1.2)
$$
compatible with (3.1.1), because
$ X_{\bC} \,(:= X\otimes_{K}\bC) $ is identified with the closed
fiber of
$ X_{R'}\otimes_{k_{R'}}\bC $ over the closed point of
$ S'_{\bC} $ defined by the inclusion
$ R' \subset \bC $.
So we get the canonical functor
$$
\iota : \cM(X/K)_{\ak} \to \MHM(X_{\bC})
\leqno(3.1.3)
$$
compatible with (3.1.1), (3.1.2).
(We will omit
$ /K $ if
$ K = \bC $.)

We denote by
$ H^{i} $ the (shifted) cohomology functor from the derived categories
of
$ \cM(X_{R'}/k_{R'}) $ and
$ \cM(X/K)_{\ak} $ to
$ \cM(X_{R'}/k_{R'})[-d_{R'}] $ and
$ \cM(X/K)_{\ak} $ respectively.
It is compatible with (3.1.1--3).

In the case
$ X = \bC $,
we define
$$
\cM_{K,\ak} = \cM(\Spec K/K)_{\ak},
$$
and it will be denoted by
$ \cM_{\ak} $ if
$ K = \bC $.
In the case
$ \cM(X) = \MHM(X) $,
it is denoted by
$ \MHS_{K,\ak} $,
and by
$ \MHS_{\ak} $ when
$ K = \bC $.

We have the canonical functor
$$
\iota : \cM_{K,\ak} \to \MHS,
\leqno(3.1.4)
$$
where the target is the category of graded-polarizable mixed
$ \bQ $-Hodge structures in the usual sense [16].

For
$ j \in \bZ $,
we have
$ \bQ_{K,\ak}(j) \in \cM_{K,\ak} $ which is represented
by a constant variation of Hodge structure of type
$ (-j,-j) $ as usual.
This is denoted by
$ \bQ_{\ak}(j) $ if
$ K = \bC $.

\medskip\noindent
{\bf 3.2.~Remarks.} (i)
The extension groups in
$ \cM_{K,\ak} $ are too big, and the natural functor (3.1.4) is not fully
faithful, see [51], 2.5 (ii).

\medskip
(ii) Let
$ \MHS_{\nabla /\ok} $ denote the category of mixed Hodge structures
whose
$ \bC $-part is endowed with an integral connection
$ \nabla $ over
$ \ok $.
Then the functor (3.1.4) factors through
$ \MHS_{\nabla /\ok} $,
and
$ \MHS_{\ak} $ is a full subcategory of
$ \MHS_{\nabla /\ok} $.

\medskip
(iii) If
$ K $ contains
$ \ok $,
we have equivalences of categories
$$
\MHM(X/K)_{\ak} \to \MHM(X/K)_{\aok},
\quad \MHS_{K,\ak} \to \MHS_{K,\aok}.
$$
induced by
$ \MHM(X_{R'}/k_{R'}) \to \MHM(X_{R'}\otimes_{k_{R'}}\ok
/\ok) $.
See [51], 2.8.

\medskip\noindent
{\bf 3.3.~Theorem.} {\it
The category
$ \cM(X/K)_{\ak} $ is an abelian category, and there exist canonical
functors
$ f_{*} $,
$ f_{!} $,
$ f^{*} $,
$ f^{!} $,
$ \bD $,
$ \otimes $,
$ \cHom $,
etc. between the derived categories
$ D^{b}\cM(X/K)_{\ak} $ in a compatible way with the functor
$ \iota $.
}

\medskip\noindent
{\it Proof.}
This follows from (1.7) by using an analogue of the generic base change
theorem in [19].

\medskip\noindent
{\bf 3.4.~Refined cycle map.} For a
$ K $-variety
$ X $ with structure morphism
$ a_{X/K} : X \to \Spec K $,
we define
$$
\aligned
\bQ_{X/K,\ak}(j)
&= {a}_{X/K}^{*}\bQ_{K,\ak}(j)
\in D^{b}\cM(X/K)_{\ak},
\\
H^{i}(X/K,\bQ_{\ak}(j))
&= H^{i}(a_{X/K})_{*}\bQ_{X/K,\ak}(j)
\in \cM_{K,\ak}.
\endaligned
$$
The latter is represented by
$ H^{i}\pi_{*}\bQ_{X_{R}/k_{R}}(j) $ for a model
$ \pi : X_{R} \to S = \Spec R $ of
$ X $.
We define an analogue of Deligne cohomology by
$$
\aligned
{H}_{\cD}^{i}(X/K,\bQ_{\ak}(j))
&= \Ext^{i}(\bQ_{K,\ak}, (a_{X/K})_{*}\bQ_{X/K,\ak}(j))
\\
\bl(
&= \Ext^{i}(\bQ_{X/K,\ak}, \bQ_{X/K,\ak}(j))\br),
\endaligned
$$
which is isomorphic to the inductive limit of
$$
{H}_{\cD}^{i}(X_{R}/k_{R},\bQ(j)) = \Ext^{i}(\bQ_{X_{R}/k_{R}},
\bQ_{X_{R}/k_{R}}(j)).
$$
If
$ K = \bC $,
we omit
$ /K $ or
$ K $,
to simplify the notation.

If
$ X/K $ is smooth, we have the refined cycle map
$$
cl : \CH^{p}(X)_{\bQ} \to {H}_{\cD}^{2p}(X/K,\bQ_{\ak}(p))
\leqno(3.4.1)
$$
by taking the inductive limit of the cycle map in (2.1.2) for models
$ X_{R}/k_{R} $ of
$ X $:
$$
cl_{R} : \CH^{p}(X_{R})_{\bQ} \to {H}_{\cD}^{2p}(X_{R}/k_{R},
\bQ_{k_{R}}(p)).
\leqno(3.4.2)
$$
This means that the cycle map is defined by taking models of cycles.

For smooth
$ K $-varieties
$ X, Y $,
let
$$
C^{i}(X,Y)_{\bQ} = \oplus_{j}
\CH^{i+\dim X_{j}}(X_{j}\times_{K}Y)_{\bQ},
$$
where the
$ X_{j} $ are the irreducible components of
$ X $.
Then the cycle map induces
$$
\aligned
C^{i}(X,Y)_{\bQ}
&\to \Ext^{2i+2\dim X}(\bQ_{X\times Y/K,\ak},
\bQ_{X\times Y/K,\ak}(i+\dim X))
\\
&=\Hom((a_{X/K})_{*}\bQ_{X/K,\ak}, (a_{Y/K})_{*}\bQ_{Y/K,\ak}(i)[2i]),
\endaligned
\leqno(3.4.3)
$$
where
$ X $ is assumed connected.
This is compatible with the composition of correspondences,
see [45], II, 3.3.
In particular, the action of
$ C^{i}(X,Y)_{\bQ} $ on the Chow groups corresponds by the cycle map
to the composition of morphisms with the image of (3.4.3), i.e. for
$ \Gamma \in C^{i}(X,Y)_{\bQ} $,
we have the commutative diagram:
$$
\CD
\CH^{p}(X)_{\bQ} @>{cl_{*}}>>
\Hom(\bQ_{K,\ak},(a_{X/K})_{*}\bQ_{X/K,\ak}(p)[2p])
\\
@VV{\Gamma_{*}}V @VV{\Gamma_{*}}V
\\
\CH^{p+i}(Y)_{\bQ} @>{cl_{*}}>> \Hom(\bQ_{K,\ak},
(a_{Y/K})_{*}\bQ_{Y/K,\ak}(p+i)[2p+2i])
\endCD
\leqno(3.4.4)
$$

We have similarly the refined cycle map for the higher Chow groups
$$
cl : \CH^{p}(X,m)_{\bQ} \to {H}_{\cD}^{2p-m}(X/K,\bQ_{\ak}(p))
\leqno(3.4.5)
$$
as the limit of (2.1.4).

\medskip\noindent
{\bf 3.5.~Indecomposable higher Chow groups.} Let
$ m = 1 $.
Then an element of
$ \CH^{p}(X,1) $ is represented by
$ \sum_{i} (Z_{i},g_{i}) $ where the
$ Z_{i} $ are integral closed subvarieties of codimension
$ p -1 $ in
$ X $ and the
$ g_{i} $ are nonzero rational functions on
$ Z_{i} $ such that
$ \sum_{i} \div g_{i} = 0 $ in
$ X $.
In particular, we have a well-defined morphism
$$
\CH^{p-1}(X)_{\bQ}\otimes_{\bZ}k^{*} \to
\CH^{p}(X,1)_{\bQ}.
\leqno(3.5.1)
$$
Its image and cokernel are denoted by
$ \CH_{\dec}^{p}(X,1)_{\bQ} $ and
$ \CH_{\ind}^{p}(X,1)_{\bQ} $ respectively.
Their elements are called decomposable and indecomposable higher cycles
respectively.

\medskip\noindent
{\bf 3.6.~Leray Filtration.} We have the Leray spectral sequence
$$
{E}_{2}^{i,j} = \Ext^{i}\bl(\bQ_{k}, H^{j}(X/K,\bQ_{\ak}(p))\br)
\Rightarrow {H}_{\cD}^{i+j}(X/K,\bQ_{\ak}(p)),
$$
which degenerates at
$ E_{2} $.
Indeed, by the decomposition theorem (see e.g. [44], 4.5.3),
we have a noncanonical isomorphism
$$
(a_{X/K})_{*}{a}_{X/K}^{*}\bQ_{K,\ak} \simeq \oplus_{j}
H^{j}(X/K,\bQ_{\ak})[-j]\quad \text{in }D^{b}\cM_{K,\ak}.
\leqno(3.6.1)
$$
We denote by
$ F_{L} $ the associated filtration on
$ {H}_{\cD}^{i+j}(X/K,\bQ_{\ak}(p)), $ and also the filtration on
$ \CH^{p}(X)_{\bQ} $ induced by the cycle map (3.4.1).
This means that
$ F_{L}^{r+1}\CH^{p}(X)_{\bQ} $ is the kernel of
$$
cl : F_{L}^{r}\CH^{p}(X)_{\bQ}
\to \Ext^{r}\bl(\bQ_{K,\ak}, H^{2p-r}(X/K,\bQ_{\ak}(p)\br),
$$
and the cycle map induces injective morphisms
$$
\Gr_{F_{L}}^{r}cl : \Gr_{F_{L}}^{r}\CH^{p}(X)_{\bQ} \to
\Ext^{r}\bl(\bQ_{K,\ak}, H^{2p-r}(X/K,\bQ_{\ak}(p))\br).
$$

\medskip\noindent
{\bf 3.7.~Remark.} By definition,
$ F_{L}^{1}\CH^{p}(X)_{\bQ} $ coincides with the subgroup
$ \CH_{\hom}^{p}(X)_{\bQ} $ consisting of cohomologically equivalent to
zero cycles.
For
$ p = 2 $,
let
$ \CH_{\AJ}^{p}(X)_{\bQ} $ denote the kernel of the Abel-Jacobi
map.
Then
$$
F_{L}^{2}\CH^{p}(X)_{\bQ} \subset \CH_{\AJ}^{p}(X)_{\bQ}.
\leqno(3.7.1)
$$
Indeed, for a model
$ \pi : X_{R} \to S = \Spec R $ of
$ X $,
we have a commutative diagram
$$
\CD
F_{L}^{1}\CH^{p}(X_{R})_{\bQ} @>>> \Ext^{1}\bl(\bQ_{S,\ak},
H^{2p-1}\pi_{*}\bQ_{X_{R}/k_{R}}(p))\br)
\\
@VVV @VVV
\\
F_{L}^{1}\CH^{p}(X_{\bC})_{\bQ} @>>>
\Ext^{1}\bl(\bQ, H^{2p-1}(X_{\bC},\bQ(p))\br),
\endCD
$$
where the vertical morphisms are induced by
$ X_{\bC} \to X_{R} $ which is given by the inclusion
$ R \to \bC $.
For
$ p = \dim X $, the equality holds in (3.7.1).
This follows from Murre's Chow-K\"unneth decomposition [36]
by using the compatibility of the action of a correspondence (3.4.4),
see e.g. [51], 3.6.
This can be generalized to algebraically equivalent to zero cycles of
any codimension, see (3.9).

We have
$ \Gr_{F_{L}}^{r}\CH^{p}(X)_{\bQ} = 0 $ for
$ r > p $ using (3.4.4), see [45], II.
It seems that the filtration
$ F_{L} $ coincides with a new filtration of M.~Green which was
explained in the conference [28].
It is expected that
$ F_{L} $ gives a conjectural filtration of Beilinson [3] and Bloch
(see also [7], [31]).

\medskip\noindent
{\bf 3.8.~Murre's filtration.} Assume that a smooth projective variety
$ X $ admits the Chow-K\"unneth decomposition in the sense of Murre
[36], [37].
Then
$ \CH^{p}(X)_{\bQ} $ has Murre's filtration
$ F_{M} $.
(See also [31].)
For
$ F_{L} $ as in (3.6), we can show (see [51], 4.9):
$$
F_{M} \subset F_{L}\quad \text{and}\quad F_{M} = F_{L}\,\,\,
\text{mod} \cap_{i} F_{L}^{i}.
\leqno(3.8.1)
$$
In particular,
$ F_{M} = F_{L} $ if the cycle map (3.4.1) for
$ X $ is injective.
The injectivity of (3.4.1) can be used for the construction
of the Chow-K\"unneth decomposition.

In [52], Shuji Saito has constructed a filtration
$ F_{\text{\rm Sh}} $ on
$ \CH^{p}(X)_{\bQ} $ by induction.
If we modify slightly his definition or assume the standard conjectures,
we can show that his filtration is contained in
$ F_{L} $ and they coincide in the case where the K\"unneth components
of the diagonal are algebraic and the cycle map (3.4.1) is injective for
the given
$ X $ (see [51], 4.9, and also [52, II] where we assume the last two
hypotheses for any smooth projective varieties).

\medskip\noindent
{\bf 3.9.~Proposition.} {\it
Let
$ \CH_{\alg}^{p}(X)_{\bQ} $ denote the subgroup consisting of
algebraically equivalent to zero cycles.
Then we have}
$$
F_{L}^{2}\CH_{\alg}^{p}(X)_{\bQ} =
\CH_{\AJ}^{p}(X)_{\bQ} \cap
\CH_{\alg}^{p}(X)_{\bQ}.
\leqno(3.9.1)
$$

\noindent
{\it Proof.}
Let
$ J^{p}(X)_{\alg} $ be the image of
$ \CH_{\alg}^{p}(X) $ by the Abel-Jacobi map to
$ J^{p}(X) $.
Then it has a structure of abelian variety.
Furthermore, if
$ X $ and a cycle
$ \zeta $ are defined over a subfield
$ K $ of
$ \bC $ which is finitely generated over
$ k $,
then so are
$ J^{p}(X)_{\alg} $ and the image of
$ \zeta $ by the Abel-Jacobi map (enlarging
$ K $ if necessary), see (3.10) below.
Since the functor (3.1.3) is induced by the inclusion
$ R' \to \bC $ which gives a geometric generic point of
$ \Spec R' $,
we get the assertion.

\medskip\noindent
{\bf 3.10.~Algebraic part of the intermediate Jacobian.} Let
$ Y $ be a closed subvariety of
$ X $ with pure codimension
$ p - 1 $,
and
$ \tY \to Y $ a resolution of singularities.
We have the Gysin morphism
$$
H^{1}(\tY,\bZ)(1-p) \to H^{2p-1}(X,\bZ).
\leqno(3.10.1)
$$
We assume that
$ Y $ is sufficiently large so that the image is maximal dimensional.
Then
$ J^{p}(X)_{\alg} $ coincides with the image of the induced morphism (see
(2.2.1)):
$$
J^{1}(\tY) \to J^{p}(X).
\leqno(3.10.2)
$$
Indeed, for
$ \zeta \in \CH_{\alg}^{p}(X) $,
there is a one-parameter family
$ \{\zeta_{s}\} $ over a connected curve
$ S $ such that
$ \zeta $ is the difference of
$ \zeta_{s_{1}} $ and
$ \zeta_{s_{0}} $ for
$ s_{0}, s_{1} \in S $,
and we may assume that
$ Y $ contains the supports of
$ \zeta_{s} $,
because the images of (3.10.1-2) do not change by enlarging
$ Y $ by the above assumption on
$ Y $.
More precisely,
$ \{\zeta_{s}\} $ comes from a cycle
$ \Gamma $ on
$ S \times X $,
and
$ \zeta $ belongs to the image of the composition of a correspondence
$ \Gamma' \in \CH^{1}(S \times \tY) $ and the pushforward by
$ \tY \to X $,
where
$ \Gamma' $ is obtained by pulling back
$ \Gamma \in \CH^{1}(S \times Y) $.

In particular,
$ J^{p}(X)_{\alg} $ has a structure of abelian variety as a quotient of
$ \Pic(\tY)^{0} $.
(See also [2]).
If
$ X $ and
$ \zeta $ are defined over
$ K $,
we may assume that so are
$ Y $ and
$ \tY $ by enlarging
$ K $ if necessary.
Then the quotient of
$ \Pic(\tY)^{0} $ is also defined over
$ K $.
Indeed, the quotient corresponds to a quotient system of realizations of
$ H^{1}(\tY_{K}/K,\bZ) $,
which is defined by using the image of (3.10.1).
Here
$ \tY_{K} $ is a model of
$ \tY $ over
$ K $,
and
$ H^{1}(\tY_{K}/K,\bZ) $ is a system of realizations with integral
coefficients, which is defined in a similar way to (1.6) and (1.8).
So the assertion follows because the divisor class is defined in
$ \Pic(\tY_{K})^{0} $.

\medskip\noindent
{\bf 3.11.~Universal ind-abelian quotient.}
We can give a more precise description of
$ \CH_{\alg}^{p}(X) $ related to [35], [42], when the base field
$ k $ is an algebraically closed field of characteristic
$ 0 $.
By the above argument, we have
$$
\CH_{\alg}^{p}(X) = \varinjlim \CH_{\alg}^{1}(Y),
\leqno(3.11.1)
$$
where the inductive limit is taken over closed subvarieties
$ Y $ of pure codimension
$ p - 1 $,
and
$ \CH_{\alg}^{1}(Y) $ is the image of
$ \CH_{\alg}^{1}(\tY) = \Pic(\tY)^{0} $ with
$ \pi : \tY \to Y $ a resolution of singularities.
As is well-known,
$ \Pic(\tY)^{0} $ is the group of
$ k $-valued points of the Picard variety.
The latter will be denoted by
$ P_{Y} $,
because it is independent of the choice of
$ \tY $.
It is the product of
$ P_{Y_{i}} $ for the irreducible components
$ Y_{i} $ of
$ Y $.

Let
$ \Lambda $ be the set of closed subvarieties
$ Y $ of pure codimension
$ p - 1 $ in
$ X $.
It has a natural ordering by the inclusion relation.
For
$ Y, Y' \in \Lambda $ such that
$ Y < Y' $,
there is a natural morphism of abelian varieties
$ \lambda _{Y,Y'} : P_{Y} \to P_{Y'} $.
So we get an inductive system of abelian varieties
$ \{P_{Y}\} $.
We say that an inductive subsystem of closed group subschemes
$ \{P'_{Y}\} $ of
$ \{P_{Y}\} $ is strict if
$ {\lambda }_{Y,Y'}^{-1}(P'_{Y'}) = P'_{Y} $ for any
$ Y < Y' $.
For example, a regular morphism to an abelian variety [42] defines
a strict subsystem by taking the kernel.
Let
$ \{P'_{Y}\} $ be the minimal strict inductive subsystem of closed
group subschemes
$ \{P_{Y}\} $ such that
$ P'_{Y}(k) $ contains the kernel of
$ P_{Y}(k) \to \CH_{\alg}^{p}(X) $ for any
$ Y $.
Let
$ \oP_{Y} = P_{Y}/P'_{Y} $.
Then
$ \{\oP_{Y} \} $ is an inductive system with injective transition
morphisms.
Let
$$
A^{p}(X)^{\ab} = \varinjlim \oP_{Y}(k).
\leqno(3.11.2)
$$
This is called the universal ind-abelian quotient of
$ \CH_{\alg}^{p}(X) $.
Clearly the ind-object
$ \{\oP_{Y}\} $ has the universal property for the regular morphisms
to abelian varieties (see [42]).
By Murre [35],
$ A^{p}(X)^{\ab} $ is an abelian variety if
$ p = 2 $.
(Indeed,
$ \dim \overline{P}_{Y} $ is bounded by using (19) of loc. cit.
and taking an abelian subvariety of
$ P_{Y} $ whose intersection with
$ P'_{Y} $ is finite.)
It coincides with the algebraic part of the intermediate Jacobian
if furthermore
$ k = \bC $.

\bigskip\bigskip
\centerline{{\bf 4. Main Results}}

\bigskip\noindent
In this section
$ X $ is assumed to be a smooth complex projective variety.
The first result is the nontriviality of the refined cycle map restricted
to the kernel of the Abel-Jacobi map, which is inspired by
Voisin's result [58] (see also [1]):

\medskip\noindent
{\bf 4.1.~Theorem.} {\it
The image of
$ \Gr_{F_{L}}^{2}cl $ in (3.6) for
$ p = \dim X $ is an infinite dimensional
$ \bQ $-vector space if
$ \Gamma (X,\Omega_{X}^{2}) \ne 0 $.
(See [51], 4.4.)
}

\medskip
This follows from Bloch's diagonal argument combined with Murre's
Chow-K\"unneth decomposition for the Albanese motive.
If
$ \Gamma (X,\Omega_{X}^{1}) = 0 $,
the latter is not necessary, and the argument is rather simple.
It is remarked by the referee that a similar assertion holds for the
cycle map to the arithmetic de Rham cohomology
$ H_{\DR}^{2p}(X/k) $ of
$ X $ which is isomorphic to the inductive limit of the de Rham
cohomology
$ H_{\DR}^{2p}(X_{R}/k) $ over the models
$ X_{R}/R $ of
$ X/\bC $.

We say that a smooth (resp. smooth projective)
$ k $-variety
$ Y $ is a {\it
$ k $-smooth} (resp. {\it
$ k $-smooth projective) model} of a complex
algebraic variety
$ X $,
if
$ Y $ has a morphism to an integral
$ k $-variety whose geometric generic fiber over
$ \Spec \bC $ is isomorphic to
$ X $.
The main point in the theory of arithmetic mixed sheaves is that the
injectivity of the refined cycle map can be reduced to that of the
Abel-Jacobi map for smooth projective models over number fields (2.1.3).
(This shows that an additional assumption in [1] is unnecessary.)

\medskip\noindent
{\bf 4.2.~Theorem.} {\it
Assume
$ k $ is a number field.
Then the cycle map
$ \CH^{p}(X)_{\bQ} \to {H}_{\cD}^{2p}(X,\bQ_{\ak}(p)) $ in (3.4.1) is
injective if the above Abel-Jacobi map (2.1.3) over
$ k_{Y} $ is injective for any
$ k $-smooth models
$ Y $ of
$ X $,
where
$ k_{Y} $ is the algebraic closure of
$ k $ in
$ \Gamma (Y,\cO_{Y}) $.
In the case
$ p = 2, $ the last assumption is reduced to the same injectivity for
any
$ k $-smooth projective models
$ Y $ of
$ X $.
(See [51], 4.6.)
}

\medskip
We can describe the image of the usual cycle map to Deligne cohomology
assuming the Hodge conjecture, if
$ H^{2p-1}(X,\bQ_{\ak}) $ is {\it global section-free} in
the sense that the local system on
$ S_{\bC} $ underlying the representative of
$ H^{2p-1}(X,\bQ_{\ak}) $ on any
$ S = \Spec R $ has no nontrivial global sections.

\medskip\noindent
{\bf 4.3.~Proposition.} {\it
If
$ H^{2p-1}(X,\bQ_{\ak}) $ is global section-free in the above sense, and
the Hodge conjecture for codimension
$ p $ cycles holds for any
$ k $-smooth projective models of
$ X $,
then the image of the usual cycle map to Deligne cohomology
$$
cl : \CH^{p}(X)_{\bQ} \to {H}_{\cD}^{2p}(X,\bQ(p))
\leqno(4.3.1)
$$
coincides with
$ \Im\bl({H}_{\cD}^{2p}(X,\bQ_{\ak}(p)) \to
{H}_{\cD}^{2p}(X,\bQ(p))\br) $,
and similarly for the Abel-Jacobi map with Deligne cohomology replaced by
the intermediate Jacobian.
(See [51], 4.1.)
}

\medskip
We can show some assertions for higher Chow groups corresponding to
(4.1--2).
Recall that the coniveau filtration
$ N^{p}H^{i}(X,\bQ) $ for a smooth proper variety
$ X $ is defined to be the kernel of
$ H^{i}(X,\bQ) \to H^{i}(U,\bQ) $ for a sufficiently
small open subvariety
$ U $ of
$ X $ such that
$ \dim X \setminus U \le p $.

\medskip\noindent
{\bf 4.4.~Theorem.} {\it
Assume
$ N^{p-2}H^{2p-3}(X,\bQ) \ne 0 $.
Then for any positive integer
$ m $,
there exist
$ \zeta_{i} \in \CH_{\hom}^{p-1}(X)_{\bQ} $ for
$ 1 \le i \le m $ together with a finitely generated subfield
$ K $ of
$ \bC $ such that for any complex numbers
$ \alpha_{1}, \dots, \alpha_{m} $ not algebraic over
$ K $,
the images of
$ \zeta_{i}\otimes \alpha_{i} $ for
$ 1 \le i \le m $ by the composition of (3.5.1) and (3.4.5) are linearly
independent over
$ \bQ $.
See [51], 5.2.
}

\medskip\noindent
{\bf 4.5.~Theorem.} {\it
Assume
$ k $ is a number field.
Then the cycle map (3.4.5) for
$ p = 2, m = 1 $ is injective if the generalized Abel-Jacobi map (2.1.5)
for the same
$ p, m $ over
$ k_{Y} $ is injective for any
$ k $-smooth projective models
$ Y $ of
$ X $,
where
$ k_{Y} $ is as in (4.2).
(See [51], 5.3.)
}

\medskip
Related to the countability of the indecomposable higher Chow group
$ \CH_{\ind}^{2}(X,1)_{\bQ} $ in (3.5), we have

\medskip\noindent
{\bf 4.6.~Theorem.} {\it
If the cycle map (3.4.5) for
$ X $ is injective for
$ p = 2, m = 1 $,
then it induces an injective morphism (see [51], 5.7):}
$$
\CH_{\ind}^{2}(X,1)_{\bQ} \to
\Ext^{1}\bl(\bQ_{\ak},H^{2}(X,\bQ_{\ak}(2))/
N^{1}H^{2}(X,\bQ_{\ak}(2))\br).
\leqno(4.6.1)
$$

For the image of (4.6.1), we can show

\medskip\noindent
{\bf 4.7.~Proposition.} {\it
Assume
$ k $ is a number field.
Then (4.6.1) is surjective if
$ H^{2}(X,\bQ_{\ak})/ $
$ N^{1}H^{2}(X,\bQ_{\ak}) $ is global section-free as in (4.3)
and if the Abel-Jacobi map (2.1.3) is injective for codimension
$ 2 $ cycles on any
$ k $-smooth projective models of
$ X $.
(See [51], 5.11.)
}

\medskip\noindent
{\bf 4.8.~Remark.} By (4.6), Voisin's conjecture on the countability of
$ \CH_{\ind}^{2}(X,1)_{\bQ} $ (see [57]) can be reduced to the
injectivity of (2.1.5) (i.e. to the hypothesis of (4.5)),
because we have an analogue of the rigidity argument due to Beilinson [3]
and M\"uller-Stach [34] as follows:

\medskip\noindent
{\bf 4.9.~Proposition.} {\it
For a smooth complex projective variety
$ X $,
let
$ \tN^{i}H^{2i}(X, \bQ_{\ak}) $ be the maximal subobject
of
$ H^{2i}(X, \bQ_{\ak}) $ which is isomorphic to a direct sum of
copies of
$ \bQ_{\ak}(-i) $.
Then the image of the morphism
$$
\CH^{p}(X,1)_{\bQ} \to \Ext^{1}(\bQ_{\ak},
H^{2p-2}(X, \bQ_{\ak}(p))/\tN^{p-1}H^{2p-2}(X, \bQ_{\ak}(p)))
\leqno(4.9.1)
$$
induced by the cycle map is countable.
(See [51], 5.9.)
}

\medskip\noindent
{\bf 4.10.~Remark.} We have the reduced higher Abel-Jacobi map
$$
\CH_{\ind}^{p}(X,1)_{\bQ} \to
\Ext^{1}\bl(\bQ, (H^{2p-2}(X,\bQ)/\Hdg^{p-1})(p)\br)
\leqno(4.10.1)
$$
in the usual sense, where
$ \Hdg^{p-1} $ denotes the group of Hodge cycles.
This is an analogue of (4.6.1).
By A. Beilinson [3] and M. Levine [33], it can be described quite
explicitly by using currents in a similar way to Griffiths' Abel-Jacobi
map (see also [24], [34]).
This is generalized to the nonproper smooth case [50].
It is not easy to construct nontrivial indecomposable higher cycles
(see [13], [14], [15], [24], [34], [57], etc.) nor nontrivial elements
in the image of (4.6.1).
An example of a nontrivial higher cycle whose support is the one point
compactification of
$ \bC^{*} $ is given in [50].
This is also an example such that its image by (4.10.1) is not contained
in the image of
$ F^{1}H^{2}(X,\bC) $ (see [15] for another example).
It was originally considered in order to find an indecomposable higher
cycle on a self-product of an elliptic curve of CM type, where the cycle
map to real Deligne cohomology used in [24] does not work.

By [38], [40], the kernel of (4.10.1) is isomorphic to
$$
\Coker(K_{2}(\bC(X))_{\bQ} \to \varinjlim
\Hom_{\MHS}(\bQ, H^{2}(U,\bQ)(2))),
\leqno(4.10.2)
$$
where the morphism is given by
$ d \log \wedge d \log $ at the level of integral logarithmic forms, and
the inductive limit is taken over nonempty open subvarieties
$ U $ of
$ X $.
This isomorphism follows easily from the localization sequence of
mixed Hodge structures together with the fact that the residue of
$ d \log f \wedge d \log g $ coincides with the differential of the
tame symbol of
$ \{f,g\} $ up to sign.
It is conjectured by Beilinson that (4.10.2) should vanish.

\bigskip\bigskip
\centerline{{\bf 5. Examples}}

\bigskip\noindent
By (4.2) and (4.5), some major problems are reduced to the injectivity of
the generalized Abel-Jacobi maps (2.1.3) and (2.1.5), and this is the
crucial point in the theory of arithmetic mixed sheaves.
However, it is not easy to verify, for example, the injectivity of
(2.1.3) even for a surface
$ X $ unless
$ p_{g}(X) = 0 $.
The problem seems to be of arithmetic nature, because this injectivity
does not hold unless the ground field
$ k $ is a number field.
We try to illustrate the difficulty of the problem in the following
example.

\medskip\noindent
{\bf 5.1.~Product of elliptic curves.} Let
$ E_{i} $ be an elliptic curve with the origin
$ O_{i} $ defined over a number field
$ k $ for
$ i = 1, 2 $.
Let
$ X = E_{1} \times E_{2} $ where the subscript
$ k $ is omitted to simplify the notation because the base change by
$ k \to \bC $ is not used here.
The choice of the origin gives a double cover
$ E_{i} \to \bP^{1} $ ramified over four points
$ \Sigma_{i} \subset \bP^{1} $,
which contain the image of the origin.

Let
$ P_{i} $ be
$ k $-valued point of
$ E_{i} $,
and
$ \zeta_{i} = [P_{i}] - [-P_{i}] $ for
$ i = 1, 2 $.
Put
$ \zeta = \zeta_{1} \times \zeta_{2} $.
Then it belongs to the kernel of the Abel-Jacobi map (i.e. of the Albanese
map in this case).
So we have to show that a multiple of
$ \zeta $ is rationally equivalent to zero.
Consider the involution
$ \sigma $ of
$ X $ defined by
$ x \to -x $.
Then
$ \zeta $ is invariant by this involution, and is identified with a cycle
of
$ X' = X/\sigma $.
This
$ X' $ is a double covering of
$ S := \bP^{1} \times \bP^{1} $ ramified over the divisor
$ D = \Sigma_{1} \times \bP^{1} \cup \bP^{1} \times
\Sigma_{2} $,
and has 16 ordinary double points as well-known.
Let
$ \tau $ be the involution associated with the double covering.
Then
$ \zeta $ is a
$ \tau $-anti-invariant cycle, i.e.
$ \tau^{*}\zeta = -\zeta $.
(Note that any cycle with rational coefficients on
$ X' $ coincides with a
$ \tau $-anti-invariant cycle modulo
$ \tau $-invariant cycles, and
$ \tau $-invariant cycles are essentially trivial modulo rational
equivalence.)

If we consider a curve which is invariant by
$ \tau $,
the description of a rational function on it is complicated.
So we try to find a curve
$ C $ on
$ S $ together with a rational function
$ g $ on
$ S $ such that
$ C $ can be lifted to a curve
$ C' $ on
$ X' $ which is birational to
$ C $,
and the divisor of the pull-back
$ g' $ of
$ g $ to
$ C' $ coincides with a multiple of
$ \zeta $ in
$ X' $.
Since
$ S $ is a self-product of
$ \bP^{1} $,
a curve on
$ S $ is described explicitly by using an equation.
However, it is not easy to express the condition on birational lifting
of the curve.
Let
$ x, y $ be affine coordinates of
$ \bA^{2} $,
which is the complement of two irreducible components of
$ D $ in
$ \bP^{1}\times \bP^{1} $.
For
$ i = 1, 2 $, let
$ f_{i} $ be the defining equation of
$ \Sigma_{i} \cap \bA^{1} $ which is a polynomial of degree
$ 3 $.
Then the restriction of
$ X' $ over
$ \bA^{2} $ is given by
$ z^{2} = f_{1}(x)f_{2}(y) $,
and an example of a liftable curve by
$ g_{1}(x,y)^{2} = f_{1}(x)f_{2}(y)g_{2}(x,y)^{2} $,
where
$ g_{1}(x,y), g_{2}(x,y) $ are polynomials with no common factors.
However, how to choose a rational function on
$ C $ is still a nontrivial problem.
If we take
$ x - c $ for
$ c \in k $ such that $ \{x = c \} $ is not contained in
$ \Sigma_{1} $, then we see that
$ [Q] - [O] $ and hence
$ [Q] - [\tau Q] $ are rationally equivalent to zero, where
$ Q $ is the sum of the points (counted with multiplicity) in the
intersection of
$ C' $ with the elliptic curve
$ \{c\}\times E_{2} $ (the sum is taken in the elliptic curve
$ E_{2} $), and
$ O $ is any of the points in the inverse image of
$ D $ which are rationally equivalent to each other.
This
$ Q $ depends only on the bidegrees of
$ g_{1}, g_{2} $, because the parameter spaces are
rational.
In particular, we get only a countable number of
$ Q $ for each
$ c $.
It is unclear whether
$ Q $ is a nontorsion point of
$ E_{2} $.

\medskip\noindent
{\bf 5.2.~Higher cycles on an elliptic curve.} Since the higher extension
groups should vanish in the category of (conjectural) mixed motives over
a number field, it is expected that a higher cycle of the form
$ \zeta \otimes \alpha $ for
$ \zeta \in \CH_{\hom}^{1}(X) $ and
$ \alpha \in k^{*} $ (see (3.5)) vanishes in
$ \CH^{2}(X,1)_{\bQ} $ if
$ X $ is smooth proper over a number field
$ k $ (see [8], [39]).

In the case of an elliptic curve
$ E $ with origin
$ O $ defined over
$ k $, take
$ P \in E(k) $, and
let
$ P_{m} = mP \in E(k) $ for
$ m = 0, 1, 2, 3 $.
Then we have a rational function
$ f $ on
$ E $ such that
$ \div f = 2[P_{1}] - [P_{0}] - [P_{2}] $.
Let
$ T_{Q} $ denote the translation by
$ Q \in E(k) $,
and
$ g = {T}_{-P}^{*}f $.
Considering the tame symbol of
$ \{cf^{3},cg^{3}\} $ for an appropriate
$ c \in k^{*} $, we see that
$ ([P_{3}] - [P_{0}])\otimes \alpha $ vanishes in
$ \CH^{2}(X,1) $ for some
$ \alpha \in k^{*} $.
We can verify that
$ \alpha \ne 1 $ for a general
$ P $ as follows.

We take the Weierstrass equation
$ y^{2} = x^{3} + Ax + B $ so that the origin of
$ E $ is the point at infinity.
Then
$ f $ is given by the function
$ T_{-P}^{*}(x - a) $ where
$ a = x(P) $ (the value of
$ x $ at
$ P) $.
Consider the function
$ (x - a)^{2}T_{-P}^{*}x $.
This can be extended to a function on a neighborhood of
$ P $.
We denote its value at
$ P $ by
$ h(P) $.
Then
$ h(P) = h(-P) $ because
$ x(P) = x(-P) $ and
$ h(P) = \lim_{Q\to P}(x(Q)-x(P))^{2}x(Q-P) $.
Furthermore, the above
$ c $ is given by
$ h(P)^{-1} $,
and
$ \alpha $ by
$ x(2P)^{9}h(P)^{-3} $.
It is easy to see that the last function of
$ P $ goes to the infinity when
$ P $ approaches to a point
$ P' $ such that
$ 2P' = O $.
(As remarked by the referee, the above
$ h(P) $ coincides with
$ y(P)^{2} $ by using the Weierstrass
$ \frak p $-function.
He also notes that we get a similar result by calculating simply the
tame symbol of
$ \{f,g\} $ and using the well-definedness of
$ \CH^{1}(X)\otimes\bC \to \CH^{2}(X,1) $,
because the tame symbol is bilinear.)

\medskip
With the notation of (5.1), a cycle
$ \zeta $ of the form
$ \zeta _{1}\times \zeta _{2} $ with
$ \zeta _{i} \in \CH^{1}(E_{i}) $ is called decomposable.
Assume that the
$ \zeta _{i} $ are homologically equivalent to zero.
Then
$ \zeta $ should be rationally equivalent to zero if the ground field
$ k $ is a number field.
But the situation is rather complicated in general.
We explain here a special case where the cycle is detected by the
refined cycle map.

\medskip\noindent
{\bf 5.3.~Strictly decomposable cycles.} Let
$ X_{1} $ and
$ X_{2} $ be smooth complex projective varieties defined over a subfield
$ k $ of
$ \bC $, i.e. there exist smooth projective
$ k $-varieties
$ X_{i,k} $ such that
$ X_{i} = X_{i,k}\otimes_{k}\bC $.
Set
$ X = X_{1} \times X_{2}, X_{k} = X_{1,k} \times_{k} X_{2,k} $.
We say that a cycle
$ \zeta $ on
$ X $ is strictly decomposable if there exist subfields
$ K_{i} $ of
$ \bC $ finitely generated over
$ k $, together with cycles
$ \zeta_{i} $ on
$ X_{i,K_{i}} := X_{i,k} \otimes_{k} K_{i} $ for
$ i = 1, 2 $ such that the algebraic closure
$ k' $ of
$ k $ in
$ K_{i} $ is independent of
$ i $,
the canonical morphism
$ K_{1} \otimes_{k'} K_{2} \to \bC $ is injective, and
$ \zeta $ coincides with the base change of the cycle
$ \zeta_{1} \times_{k'} \zeta_{2} $ on
$ X_{1,K_{1}} \times_{k'} X_{2,K_{2}} =
X_{k'} \otimes_{k'} (K_{1} \otimes_{k'} K_{2}) $ by
$ K_{1} \otimes_{k'} K_{2} \to \bC $.
Here we may assume
$ k' = k $, replacing
$ k $ if necessary.
We say that a strictly decomposable cycle is of bicodimension
$ (p_{1},p_{2}) $ if
$ \codim\,\zeta_{i} = p_{i} $.
Put
$ p = p_{1} + p_{2} $.

Let
$ R_{i} $ be a finitely generated smooth
$ k $-subalgebra of
$ K_{i} $ such that the fraction field is
$ K_{i} $,
and
$ \zeta_{i} $ is defined over
$ R_{i} $.
Set
$ S_{i} = \Spec R_{i}\otimes_{k}\bC $.
Let
$ \xi_{i}^{j} \in
H^{2p_{i}-j}(X_{i},\bQ)\otimes H^{j}(S_{i},\bQ)(p_{i}) $ be
the K\"unneth components of the cycle class of
$ \zeta_{i}\otimes_{k}\bC $ in
$ H^{2p_{i}}(X_{i}\times S_{i},\bQ(p_{i})) $.
Let
$ \tM_{i} $ be the pull-back of
$$
M_{i} := H^{2p_{i}-1}(X_{i},\bQ)(p_{i})
$$
by
$ a_{S_{i}} : S_{i} \to \Spec \bC $.
If
$ \xi_{i}^{0} = 0 $, then
$ \zeta_{i} \in {F}_{L}^{1}\CH^{p_{i}}(X_{i,k} \otimes_{k}
K_{i})_{\bQ} $, and
$ \Gr_{{F}_{L}}^{1}cl_{R_{i}}(\zeta_{i}) $ (see (3.4.2)) gives
$$
\widetilde{\xi}_{i}^{1} \in
\Ext^{1}(\bQ_{S_{i}},\tM_{i}).
$$
By the Leray spectral sequence, we have an exact sequence
$$
0 \to \Ext^{1}(\bQ, M_{i}) \to
\Ext^{1}(\bQ_{S_{i}},\tM_{i}) \to
\Hom(\bQ,H^{1}(S_{i},\bQ)\otimes M_{i}),
$$
where the first morphism is given by the pull-back by
$ a_{S_{i}} $.
Note that
$ \xi_{i}^{1} $ coincides with the image of
$ \widetilde{\xi}_{i}^{1} $ in the last term.
If
$ \xi_{i}^{1} = 0 $, then
$ \widetilde{\xi}_{i}^{1} $ comes from
$ \Ext^{1}(\bQ, M_{i}) $, and we may assume
$ R_{i} = k $ as long as
$ \Gr_{{F}_{L}}^{1}cl_{R_{i}}(\zeta_{i}) $ is concerned (e.g. if
$ p_{i} = 1) $.

We say that a strictly decomposable cycle
$ \zeta $ is {\it degenerate}, if either
$ \xi_{i}^{0} = \xi_{i}^{1} = 0 $ for both
$ i $, or
$ \xi_{i}^{0} = \widetilde{\xi}_{i}^{1} = 0 $ for some
$ i $.
In the case
$ p_{i} = 1 $, it is expected that a degenerate
$ \zeta $ is rationally equivalent to zero by Beilinson's conjecture
(2.3) when
$ k $ is a number field.
However, if
$ \zeta $ is nondegenerate, we can detect it by the refined cycle map
as follows.
(This is by joint work with A. Rosenschon [41]).

\medskip\noindent
{\bf 5.4.~Theorem.} {\it
With the above notation, let
$ \zeta $ be a strictly decomposable cycle of codimension
$ p $.
Then
$ cl(\zeta ) \ne 0 $ in the notation of (3.6).
More precisely,
$ \zeta \in {F}_{L}^{r}\CH^{p}(X)_{\bQ} $ and
$ \Gr_{{F}_{L}}^{r}cl(\zeta) \ne 0 $ if the number of the
$ i $ such that
$ \xi_{i}^{0} = 0 $ is
$ r $.
}

\medskip\noindent
{\bf 5.5.~Remark.} The assertion (5.4) was first considered in the case
both
$ \xi_{1}^{1} $ and
$ \xi_{2}^{1} $ are nonzero, in order to show the nonvanishing of the
composition of certain extension classes (see [49]).
A. Rosenschon studied Nori's construction of a cycle [53] in the case
of a self-product of an elliptic curve without complex multiplication,
and obtained a special case of (5.4).
Then these two were generalized to (5.4), see [41].
This can be extended to the higher cycle case (loc. cit.)
As an application, we can show that
$ \CH_{\ind}^{p+1}(X_{1}\times X_{2},1)_{\bQ} $ is uncountable if
$ \Gamma(X_{1},\Omega_{X_{1}}^{1}) \ne 0 $ and the reduced higher
Abel-Jacobi map (see (4.10.1)) for
$ X_{2} $ is not zero.
This is a generalization of a result of Gordon and Lewis [24]
(see also [14]).

\bigskip\bigskip
\centerline{{\bf References}}

\bigskip
\item{[1]}
M. Asakura, Motives and algebraic de Rham cohomology, in The arithmetic
and geometry of algebraic cycles (Banff), CRM Proc. Lect. Notes, 24, AMS,
2000, pp. 133--154.

\item{[2]}
L. Barbieri-Viale, On algebraic 1-motives related to Hodge cycles,
preprint (math.AG/ 0103179).

\item{[3]}
A. Beilinson, Higher regulators and values of
$ L $-functions, J. Soviet Math. 30 (1985), 2036--2070.

\item{[4]}
\SameAuthor, Height pairing between algebraic cycles, Lect. Notes in
Math., vol. 1289, Springer, Berlin, 1987, pp. 1--26.

\item{[5]}
\SameAuthor, On the derived category of perverse sheaves, ibid.
pp. 27--41.

\item{[6]}
A. Beilinson, J. Bernstein and P. Deligne, Faisceaux pervers,
Ast\'erisque, vol. 100, Soc. Math. France, Paris, 1982.

\item{[7]}
S. Bloch, Lectures on algebraic cycles, Duke University Mathematical
series 4, Durham, 1980.

\item{[8]}
\SameAuthor, Algebraic
$ K $-theory and classfield theory for arithmetic surfaces,
Ann. of Math. 114 (1981), 229--265.

\item{[9]}
\SameAuthor, Algebraic cycles and values of
$ L $-functions, J. Reine Angew. Math. 350 (1984), 94--108.

\item{[10]}
\SameAuthor, Algebraic cycles and higher
$ K $-theory, Advances in Math., 61
(1986), 267--304.

\item{[11]}
A. Borel et al., Algebraic D-modules, Perspectives in Math. 2, Academic
Press, 1987.

\item{[12]}
J. Carlson, Extensions of mixed Hodge structures, in Journ\'ees
de G\'eom\'etrie Alg\'ebrique d'Angers 1979, Sijthoff-Noordhoff
Alphen a/d Rijn, 1980, pp. 107--128.

\item{[13]}
A. Collino, Griffiths' infinitesimal invariant and higher
$ K $-theory on hyperelliptic Jacobians, J. Alg. Geom. 6 (1997),
393--415.

\item{[14]}
A. Collino and N. Fakhruddin, Indecomposable higher Chow cycles on
Jacobians, preprint (math.AG/0005171).

\item{[15]}
P. del Angel and S. M\"uller-Stach,
The transcendental part of the regulator map for
$ K_{1} $ on a mirror family of
$ K3 $ surfaces, preprint (math.AG/0008207).

\item{[16]}
P. Deligne, Th\'eorie de Hodge I, Actes Congr\`es Intern. Math., 1970,
vol. 1, 425-430; II, Publ. Math. IHES, 40 (1971), 5--57; III ibid., 44
(1974), 5--77.

\item{[17]}
\SameAuthor, Equations Diff\'erentielles \`a Points Singuliers
R\'eguliers, Lect. Notes in Math. vol. 163, Springer, Berlin, 1970.

\item{[18]}
\SameAuthor, Le formalisme des cycles \'evanescents, in SGA7 XIII and XIV,
Lect. Notes in Math. vol. 340, Springer, Berlin, 1973, pp. 82--115 and
116--164.

\item{[19]}
\SameAuthor, Th\'eor\`eme de finitude en cohomologie
$ l $-adique, in SGA 4 1/2, Lect. Notes in Math., vol. 569, Springer,
Berlin, 1977, pp. 233--261.

\item{[20]}
\SameAuthor, Valeurs de fonctions L et p\'eriodes d'int\'egrales, in
Proc. Symp. in pure Math., 33 (1979) part 2, pp. 313--346.

\item{[21]}
P. Deligne, J. Milne, A. Ogus and K. Shih, Hodge Cycles, Motives, and
Shimura varieties, Lect. Notes in Math., vol 900, Springer, Berlin, 1982.

\item{[22]}
F. El Zein and S. Zucker, Extendability of normal functions associated
to algebraic cycles, in Topics in transcendental algebraic geometry,
Ann. Math. Stud., 106, Princeton Univ. Press, Princeton, N.J., 1984,
pp. 269--288.

\item{[23]}
H. Esnault and E. Viehweg, Deligne-Beilinson cohomology,
in Beilinson's conjectures on Special Values of
$ L $-functions,
Academic Press, Boston, 1988, pp. 43--92.

\item{[24]}
B. Gordon and J. Lewis, Indecomposable higher Chow cycles on products
of elliptic curves, J. Alg. Geom. 8 (1999), 543--567.

\item{[25]}
M. Green, Griffiths' infinitesimal invariant and the Abel-Jacobi
map, J. Diff. Geom. 29 (1989), 545--555.

\item{[26]}
\SameAuthor, What comes after the Abel-Jacobi map? preprint.

\item{[27]}
\SameAuthor, Algebraic cycles and Hodge theory
(Lecture notes at the Banff conference, based on a collaboration with
P. Griffiths).

\item{[28]}
\SameAuthor, in this volume.

\item{[29]}
P. Griffiths, On the period of certain rational integrals I, II,
Ann. Math. 90 (1969), 460--541.

\item{[30]}
U. Jannsen, Mixed motives and algebraic
$ K $-theory, Lect. Notes in Math., vol. 1400, Springer, Berlin, 1990.

\item{[31]}
\SameAuthor, Motivic sheaves and filtrations on Chow groups, Proc. Symp.
Pure Math. 55 (1994), Part 1, pp. 245--302.

\item{[32]}
M. Kashiwara, A study of variation of mixed Hodge structure,
Publ. RIMS, Kyoto Univ. 22 (1986), 991--1024.

\item{[33]}
M. Levine, Localization on singular varieties, Inv. Math. 91 (1988),
423--464.

\item{[34]}
S. M\"uller-Stach, Constructing indecomposable motivic cohomology
classes on algebraic surfaces, J. Alg. Geom. 6 (1997), 513--543.

\item{[35]}
J.P. Murre, Applications of algebraic
$ K $-theory to the theory of algebraic cycles, Lect. Notes in Math.,
vol. 1124, Springer, Berlin, 1985, pp. 216---261.

\item{[36]}
\SameAuthor, On the motive of an algebraic surface, J. Reine Angew. Math.
409 (1990), 190--204.

\item{[37]}
\SameAuthor, On a conjectural filtration on Chow groups of an algebraic
variety, Indag. Math. 4 (1993), 177--201.

\item{[38]}
C. Pedrini, Bloch's conjecture and the
$ K $-theory of projective surfaces, in The arithmetic and geometry of
algebraic cycles, CRM Proc. Lecture Notes, 24, Amer. Math. Soc.,
Providence, 2000, pp. 195--213.

\item{[39]}
W. Raskind, On
$ K_{1} $ of curves over global fields, Math. Ann. 288 (1990),
179--193.

\item{[40]}
A. Rosenschon, Indecomposable Elements of
$ K_{1} $,
$ K $-theory 16 (1999), 185--199.

\item{[41]}
A. Rosenschon and M. Saito, Cycle map for strictly decomposable cycles,
preprint (math.AG/0104079).

\item{[42]}
H. Saito, Abelian varieties attached to cycles of intermediate dimension,
Nagoya Math. J. 75 (1979), 95--119.

\item{[43]}
M. Saito, Modules de Hodge polarisables, Publ. RIMS, Kyoto Univ. 24
(1988), 849--995.

\item{[44]}
\SameAuthor, Mixed Hodge Modules, Publ. RIMS, Kyoto Univ., 26 (1990),
221--333.

\item{[45]}
\SameAuthor, Hodge conjecture and mixed motives, I, Proc. Symp. Pure Math.
53 (1991), 283--303; II, in Lect. Notes in Math., vol. 1479, Springer,
Berlin, 1991, pp. 196--215.

\item{[46]}
\SameAuthor, Some remarks on the Hodge type conjecture, Proc. Symp. Pure
Math. 55 (1991), Part 1, pp. 85--100.

\item{[47]}
\SameAuthor, On the formalism of mixed sheaves, preprint RIMS-784, Aug.
1991.

\item{[48]}
\SameAuthor, Admissible normal functions, J. Alg. Geom. 5 (1996),
235--276.

\item{[49]}
\SameAuthor, Bloch's conjecture and Chow motives, preprint
(math.AG/0002009).

\item{[50]}
\SameAuthor, Bloch's conjecture, Deligne cohomology and higher Chow
groups, preprint RIMS-1284.

\item{[51]}
\SameAuthor, Arithmetic mixed sheaves, Inv. Math. 144 (2001), 533--569.

\item{[52]}
S. Saito, Motives and filtrations on Chow groups, Inv. Math. 125 (1996),
149 --196; II, in Proceedings
of the NATO Advanced Study Institute on The arithmetic and geometry of
algebraic cycles (B.B. Gordon et al. eds.), Kluwer Academic, Dordrecht,
2000, pp. 321-346.

\item{[53]}
C. Schoen, Zero cycles modulo rational equivalence for some varieties
over fields of transcendence degree one, Proc. Symp. Pure Math. 46 (1987),
part 2, pp. 463--473.

\item{[54]}
J. Steenbrink and S. Zucker, Variation of mixed Hodge structure I,
Inv. Math. 80 (1985), 489--542.

\item{[55]}
C. Voisin, Variations de structures de Hodge et z\'ero-cycles sur
les surfaces g\'en\'erales, Math. Ann. 299 (1994), 77--103.

\item{[56]}
\SameAuthor, Transcendental methods in the study of algebraic cycles,
in Lect. Notes in Math. vol. 1594, pp. 153--222.

\item{[57]}
\SameAuthor, Remarks on zero-cycles of self-products of varieties, in
Moduli of Vector Bundles, Lect. Notes in Pure and Applied Mathematics,
vol. 179, M. Dekker, New York, 1996, pp. 265--285.

\item{[58]}
\SameAuthor, Some results on Green's higher Abel-Jacobi map, Ann. Math.
149 (1999), 451--473.

\bigskip
\noindent
\ver

\bye